# О связи имитационной логит динамики в популяционной теории игр и метода зеркального спуска в онлайн оптимизации на примере задачи выбора кратчайшего маршрута

*Гасников Александр Владимирович*[1,2,3] gasnikov@yandex.ru
*Лагуновская Анастасия Александровна*[1,3] nastu4e4ik@yandex.ru
*Морозова Лариса Эдуардовна*[1] lmorozova@hse.ru

[1]Центр экономики транспорта и Центр исследования транспортной политики Института экономики транспорта и транспортной политики НИУ ВШЭ;
[2]ИППИ РАН; [3]ПреМоЛаб МФТИ

**Аннотация**

В работе описывается метод зеркального спуска для задач стохастической онлайн оптимизации на симплексе и прямом произведении симплексов. На базе этого метода строятся оптимальные стратегии пользователей транспортной сети при выборе маршрутов следования. Поведение всех пользователей, действующих, согласно таким стратегиям, порождает имитационную логит динамику в популяционной игре, соответствующей модели Бэкмана равновесного распределения потоков по путям. Таким образом, на конкретном примере (The Shortest Path Problem), в работе показывается связь онлайн оптимизации и популяционной теории игр. Обнаружение отмеченной связи составляет основной результат данной работы.

**Ключевые слова:** метод зеркального спуска, онлайн оптимизация, кратчайший путь в графе, равновесное распределение транспортных потоков по путям.

## 1. Введение

В литературе по онлайн оптимизации почетное место занимает, так называемая, "Задача о выборе кратчайшего пути" ("The Shortest Path Problem"), см., например, п. 5.4 [1]. Основной результат здесь заключается в описании "оптимальной" стратегии пользователя транспортной сети (на базе алгоритма "Follow the Perturbed Leader"), из дня в день выбирающего маршрут следования, исходя из истории загрузок графа транспортной сети.

В литературе по равновесной теории транспортных потоков наиболее популярными являются модели равновесного распределения потоков по путям. Одной из первых (и по-прежнему наиболее популярных) моделей такого рода является модель Бекмана [2] (также называемая BMW-моделью). Современные исследования этой модели связаны с ее пониманием, как популяционной игры загрузок (как следствие, потенциальной игры [3]), поиск равновесия (Нэша) в которой сводится к задаче выпуклой оптимизации. Упомянутый эволюционный подход, в частности, приводит к изучению различных естественных динамик (наилучших ответов, репликаторов, имитационной логит динамики и др.), отражающих "нащупывание" пользователями транспортной сети равновесия [4]. Все эти динамики положительно коррелированны с антиградиентной динамикой, поэтому все они приводят в конечном итоге к одному и тому же равновесию (или в более общем случае к одному и тому же множеству равновесий). Тем не менее, возникает желание глубже разобраться с природой этих динамик. Понять чем та или иная динамика дополнительно (помимо отражения рациональности игроков/пользователей транспортной сети) примечательна.

В данной работе мы постараемся пояснить, чем примечательна имитационная логит динамика, пояснив ее связь с алгоритмом поведения "Follow the Perturbed Leader", а точнее с переформулировкой этого алгоритма на языке современной выпуклой онлайн оптимизации: с методом зеркального спуска [5 – 11].

В разделе 2 статьи описаны различные варианты классического метода зеркального спуска применительно к задачам стохастической онлайн оптимизации с шумами. В виду последующих приложений мы сосредоточимся на множествах вида симплекса и прямого произведения симплексов. Хотя во многом эти результаты ранее были известны, тем не менее, в такой общности, в которой они приведены в данной статье, нам не удалось найти точной ссылки, поэтому было решено посвятить этому отдельный раздел 2. В разделе 3 мы используем результаты раздела 2 (в данной статье не в максимальной общности) для объяснения имитационной логит динамики, возникающей при описании поведения пользователей транспортной сети в модели Бекмана.

## 2. Метод зеркального спуска для задач стохастической онлайн оптимизации с неточным оракулом

Сформулируем основную задачу стохастической онлайн оптимизации с неточным оракулом. Требуется подобрать последовательность $\{x^k\} \in Q$ так, чтобы минимизировать псевдо регрет [6 – 11]:

$$\text{Regret}_N\left(\{f_k(\cdot)\}, \{x^k\}\right) = \frac{1}{N}\sum_{k=1}^{N} f_k(x^k) - \min_{x \in Q} \frac{1}{N}\sum_{k=1}^{N} f_k(x) \qquad (1)$$

на основе доступной информации

$$\left\{\nabla_x \tilde{f}_1(x^1, \xi^1); ...; \nabla_x \tilde{f}_{k-1}(x^{k-1}, \xi^{k-1})\right\}$$

при расчете $x^k$. Причем выполнено **условие**[1]

1. для любых $k = 1,...,N$ ($\Xi^{k-1}$ – сигма алгебра, порожденная $\xi^1, ..., \xi^{k-1}$)

$$\left\|\nabla_x \tilde{f}_k(x^k, \xi^k) - \nabla_x f_k(x^k, \xi^k)\right\|_* \le \delta,$$

$$E_{\xi_k}\left[\nabla_x f_k(x^k, \xi^k)\right] = \nabla f_k(x^k).$$

Здесь случайные величины $\{\xi^k\}$ могут считаться независимыми одинаково распределенными. Онлайновость постановки задачи допускает, что на каждом шаге $k$ функция $f_k(\cdot)$ может выбираться из рассматриваемого класса функций враждебно по отношению к используемому нами методу генерации последовательности $\{x^k\}$. В частности, $f_k(\cdot)$ может зависеть от

$$\left\{x^1, \xi^1, f_1(\cdot); ...; x^{k-1}, \xi^{k-1}, f_{k-1}(\cdot); x^k\right\}.$$

Относительно класса функций, из которого выбираются $\{f_k(\cdot)\}$, в данной работе будем предполагать выполненными следующие **условия**:

2. $\{f_k(\cdot)\}$ – выпуклые функции;
3. для любых $k = 1,...,N$, $x \in Q$

---

[1] Считаем также, что $\max_{x,y \in Q} \|x - y\| \le \tilde{R}$.

$$\left\|\nabla_x \tilde{f}_k(x,\xi)\right\|_*^2 \leq M^2.$$

Опишем метод зеркального спуска для решения задачи (1) (здесь можно следовать огромному числу литературных источников, мы в основном будем следовать работам [12, 13]). Введем норму $\|\ \|$ в прямом пространстве (сопряженную норму будем обозначать $\|\ \|_*$) и прокс-функцию $d(x)$ сильно выпуклую относительно этой нормы, с константой сильной выпуклости $\geq 1$. Выберем точку старта

$$x^1 = \arg\min_{x \in Q} d(x),$$

считаем, что

$$d(x^1) = 0,\ \nabla d(x^1) = 0.$$

Введем брэгмановское "расстояние"

$$V_x(y) = d(y) - d(x) - \langle \nabla d(x), y - x \rangle.$$

Везде в дальнейшем будем считать, что

$$d(x) = V_{x^1}(x) \leq R^2 \text{ для всех } x \in Q.$$

Определим оператор "проектирования" согласно этому расстоянию

$$\text{Mirr}_{x^k}(g) = \arg\min_{y \in Q}\left\{\langle g, y - x^k \rangle + V_{x^k}(y)\right\}.$$

Метод зеркального спуска (МЗС) для задачи (1) будет иметь вид, см., например, [13]

$$x^{k+1} = \text{Mirr}_{x^k}\left(\alpha \nabla_x \tilde{f}_k(x^k, \xi^k)\right),\ k = 1,...,N.$$

Тогда при выполнении условий (2) для любого $u \in Q$, $k = 1,...,N$ имеет место неравенство, см., например, [13]

$$\alpha \left\langle \nabla_x \tilde{f}_k(x^k, \xi^k), x_k - u \right\rangle \leq \frac{\alpha^2}{2}\left\|\nabla_x \tilde{f}_k(x^k, \xi^k)\right\|_*^2 + V_{x^k}(u) - V_{x^{k+1}}(u).$$

Это неравенство несложно получить в случае евклидовой прокс-структуры $d(x) = \|x\|_2^2/2$ [14] (в этом случае МЗС для задачи (1) есть просто вариант обычного метода проекции градиента). Разделим сначала выписанное неравенство на $\alpha$ и возьмем условное математическое ожидание $E_{\xi^k}\left[\ \cdot\ |\Xi^{k-1}\right]$, затем просуммируем то, что получится по $k = 1,...,N$, используя условие 1. Затем возьмем от того, что получилось при суммировании, полное математическое ожидание, учитывая условие 3. В итоге, выбирая $u = x_*$ (решение задачи $\sum_{k=1}^N f_k(x) \to \min_{x \in Q}$), получим при условиях 1, 2, 4 [11]

$$N \cdot E\left[\text{Regret}_N\left(\{f_k(\cdot)\},\{x^k\}\right)\right] \leq \frac{V_{x^1}(x_*)}{\alpha} - \frac{E\left[V_{x^{N+1}}(x_*)\right]}{\alpha} + \left(\frac{1}{2}M^2\alpha + \tilde{R}\delta\right)N \leq$$
$$\leq \frac{R^2}{\alpha} + \left(\frac{1}{2}M^2\alpha + \tilde{R}\delta\right)N,$$

выбирая[2]

$$\alpha = \frac{R}{M}\sqrt{\frac{2}{N}},$$

получим

---

[2] Можно получить и адаптивный вариант приводимой далее оценки, для этого потребуется использовать метод двойственных усреднений [11, 13, 14].

$$E\left[\operatorname{Regret}_N\left(\{f_k(\cdot)\},\{x^k\}\right)\right] \le MR\sqrt{\frac{2}{N}} + \tilde{R}\delta. \qquad (2)$$

Немного более аккуратные рассуждения (использующие неравенство Азума–Хефдинга) позволяют уточнить оценку (2) следующим образом (см., например, [15]):

$$\operatorname{Regret}_N\left(\{f_k(\cdot)\},\{x^k\}\right) \le M\sqrt{\frac{2}{N}}\left(R + 2\tilde{R}\sqrt{\ln(\sigma^{-1})}\right) + \tilde{R}\delta \qquad (3)$$

с вероятностью $\ge 1-\sigma$.

Оценки (2), (3) являются неулучшаемыми с точностью до мультипликативного числового множителя. Причем верно это и для детерминированных (не стохастических) постановок, в которых нет шумов ($\delta = 0$), при этом можно ограничиться классом линейных функций [1].

Рассмотрим три примера, которые понадобятся нам в дальнейшем [11, 15].

**Пример 1 (симплекс).** Предположим, что

$$Q = S_n(1) = \left\{x \ge 0: \ \sum_{i=1}^{n} x_i = 1\right\}.$$

Выберем

$$\|\ \| = \|\ \|_1, \ d(x) = \ln n + \sum_{i=1}^{n} x_i \ln x_i.$$

Тогда МЗС примет следующий вид ($\alpha = M^{-1}\sqrt{2\ln n/N}$):

$$x_i^1 = 1/n, \ i = 1,...,n,$$

при $k = 1,...,N$, $i = 1,...,n$

$$x_i^{k+1} = \frac{\exp\left(-\sum_{r=1}^{k}\alpha\frac{\partial \tilde{f}_r(x^r,\xi^r)}{\partial x_i}\right)}{\sum_{l=1}^{n}\exp\left(-\sum_{r=1}^{k}\alpha\frac{\partial \tilde{f}_r(x^r,\xi^r)}{\partial x_l}\right)} = \frac{x_i^k \exp\left(-\alpha\frac{\partial \tilde{f}_k(x^k,\xi^k)}{\partial x_i}\right)}{\sum_{l=1}^{n}x_l^k \exp\left(-\alpha\frac{\partial \tilde{f}_k(x^k,\xi^k)}{\partial x_l}\right)}.$$

Оценки псевдо регрета будут иметь вид:

$$E\left[\operatorname{Regret}_N\left(\{f_k(\cdot)\},\{x^k\}\right)\right] \le M\sqrt{\frac{2\ln n}{N}} + 2\delta,$$

$$\operatorname{Regret}_N\left(\{f_k(\cdot)\},\{x^k\}\right) \le M\sqrt{\frac{2}{N}}\left(\sqrt{\ln n} + 4\sqrt{\ln(\sigma^{-1})}\right) + 2\delta$$

с вероятностью $\ge 1-\sigma$.

**Пример 2 (прямое произведение симплексов).** Предположим, что

$$x = (z^1,...,z^m) \in Q = \prod_{j=1}^{m} S_{n_j}(d_j).$$

Выберем

$$\|x\| = \sqrt{\sum_{j=1}^{m}\|z^j\|_1^2}, \ d(x) = \sum_{j=1}^{m}d_j d^j(z^j), \ d^j(z^j) = d_j\ln n_j + \sum_{i=1}^{n_j}z_i^j\ln\left(\frac{z_i^j}{d_j}\right).$$

Тогда вводя обозначения

$$\alpha = \frac{R}{M}\sqrt{\frac{2}{N}}, \ R^2 = \sum_{j=1}^{m}d_j^2\ln n_j,$$

МЗС можно записать следующим образом:
$$z_i^j = d_j/n_j, \ i = 1,...,n_j,$$

при $k = 1,...,N$, $i = 1,...,n_j$, $j = 1,...,m$

$$z_i^{j,k+1} = d^j \frac{\exp\left(-\sum_{r=1}^k \alpha \frac{\partial \tilde{f}_r(x^r, \xi^r)}{\partial z_i^j}\right)}{\sum_{l=1}^n \exp\left(-\sum_{r=1}^k \alpha \frac{\partial \tilde{f}_r(x^r, \xi^r)}{\partial z_l^j}\right)} = d^j \frac{x_i^k \exp\left(-\alpha \frac{\partial \tilde{f}_k(x^k, \xi^k)}{\partial z_i^j}\right)}{\sum_{l=1}^n x_l^k \exp\left(-\alpha \frac{\partial \tilde{f}_k(x^k, \xi^k)}{\partial z_l^j}\right)}.$$

Оценки псевдо регрета будут иметь вид:

$$E\left[\operatorname{Regret}_N\left(\{f_k(\cdot)\}, \{x^k\}\right)\right] \leq M \sqrt{\frac{2}{N}} \sqrt{\sum_{j=1}^m d_j^2 \ln n_j} + 2\delta \sqrt{\sum_{j=1}^m d_j^2},$$

$$\operatorname{Regret}_N\left(\{f_k(\cdot)\}, \{x^k\}\right) \leq M \sqrt{\frac{2}{N}} \left(\sqrt{\sum_{j=1}^m d_j^2 \ln n_j} + 4\sqrt{\sum_{j=1}^m d_j^2} \sqrt{\ln(\sigma^{-1})}\right) + 2\delta \sqrt{\sum_{j=1}^m d_j^2}$$

с вероятностью $\geq 1 - \sigma$.

**Пример 3 (выбор среди вершин симплекса).** Вернемся к примеру 1 и будем дополнительно считать, что по условию задачи $\{x^k\}$ должны выбираться среди вершин единичного симплекса $S_n(1)$. Также как и раньше онлайновость постановки задачи допускает, что на каждом шаге $k$ функция $f_k$ может подбираться из рассматриваемого класса функций враждебно по отношению к используемому нами методу генерации последовательности $\{x^k\}$. В частности, $f_k$ может зависеть от
$$\{x^1, \xi^1, f_1(\cdot);...;x^{k-1}, \xi^{k-1}, f_{k-1}(\cdot)\},$$
и даже от распределения вероятностей $p^k$, согласно которому осуществляется выбор $x^k$. Чтобы можно было работать с таким классом задач, нам придется наложить дополнительные **условия**:

4. $f_k(x) = \langle l^k, x \rangle$, $k = 1,...,N$;

5. На каждом шаге генерирование случайной величины $x^k$ согласно распределению вероятностей $p^k$ осуществляется независимо ни от чего. Выбор $f_k$ осуществляется без знания реализации $x^k$.

Следуя примеру 1, положим $p_i^1 = x_i^1 = 1/n$, $i = 1,...,n$. При $k = 1,...,N$, $i = 1,...,n$ согласно распределению вероятностей ($\alpha = M^{-1}\sqrt{2\ln n/N}$)

$$p_i^{k+1} = \frac{\exp\left(-\sum_{r=1}^k \alpha \frac{\partial \tilde{f}_r(x^r, \xi^r)}{\partial x_i}\right)}{\sum_{l=1}^n \exp\left(-\sum_{r=1}^k \alpha \frac{\partial \tilde{f}_r(x^r, \xi^r)}{\partial x_l}\right)} = \frac{p_i^k \exp\left(-\alpha \frac{\partial \tilde{f}_k(x^k, \xi^k)}{\partial x_i}\right)}{\sum_{l=1}^n p_l^k \exp\left(-\alpha \frac{\partial \tilde{f}_k(x^k, \xi^k)}{\partial x_l}\right)},$$

генерируем случайную величину $i(k+1)$, и полагаем
$$x_{i(k+1)}^{k+1} = 1, \ x_j^{k+1} = 0, j \neq i(k+1).$$

Оценки псевдо регрета будут иметь вид:

$$E\left[\mathrm{Regret}_N\left(\{f_k(\cdot)\},\{x^k\}\right)\right] \le M\sqrt{\frac{2\ln n}{N}} + 2\delta,$$

$$\mathrm{Regret}_N\left(\{f_k(\cdot)\},\{x^k\}\right) \le M\sqrt{\frac{2}{N}}\left(\sqrt{\ln n} + 6\sqrt{\ln(\sigma^{-1})}\right) + 2\delta$$

с вероятностью $\ge 1-\sigma$.

## 3. Приложение к задаче о выборе кратчайшего пути

Рассмотрим транспортную сеть, которую будем представлять ориентированным графом $\langle V, E\rangle$, где $V$ – множество вершин, а $E$ – множество ребер. Обозначим множество пар источник-сток через $OD \subseteq V \otimes V$ ($|OD| = m$); $d_w$ – корреспонденция, отвечающая паре $w$; $x_p$ – поток по пути $p$; $P_w$ – множество путей, отвечающих корреспонденции $w$, $P = \bigcup_{w \in OD} P_w$ – множество всех путей. Обозначим через $L$ – максимальное число ребер в пути из $P$. Будем считать, что затраты на прохождения ребра $e \in E$ описываются неубывающей (и ограниченной в рассматриваемом диапазоне значений) функцией

$$0 \le \tau_e(f_e) \le \tilde{M},$$

где $f_e$ – поток по ребру $e$:

$$f_e(x) = \sum_{p \in P} \delta_{ep} x_p, \quad \delta_{ep} = \begin{cases} 1, & e \in p \\ 0, & e \notin p \end{cases}.$$

Положим $\breve{M} = \tilde{M}L$. Введем $G_p(x)$ – затраты на проезд по пути $p$:

$$G_p(x) = \sum_{e \in E} \tau_e(f_e(x))\delta_{ep}.$$

Введем также множество (прямое произведение симплексов), на котором транспортная сеть "будет жить"

$$X = \left\{x \ge 0: \sum_{p \in P_w} x_p = d_w, w \in OD\right\},$$

и функцию, порождающее потенциальное векторное поле $G(x)$:

$$\Psi(x) = \sum_{e \in E} \int_0^{f_e(x)} \tau_e(z)dz.$$

Основное свойство этой функции заключается в том, что

$$\nabla\Psi(x) = G(x).$$

Будем считать, что число пользователей транспортной сети большое:

$$d_w := d_w \cdot \bar{N}, \ \bar{N} \gg 1, \ w \in OD,$$

но в функциях затрат это учитывается

$$\tau_e(f_e) := \tau_e(f_e/\bar{N}).$$

Таким образом, далее под $d_w$, $x$, $f$ будем понимать соответствующие прошкалированные (по $\bar{N}$) величины [4].

Выберем корреспонденцию $w \in OD$, и рассмотрим пользователя транспортной сетью, соответствующего этой корреспонденции. Стратегией пользователя является выбор одного из возможных путей следования $p \in P_w$. Будем считать, что пользователь мало что знает об устройстве транспортной системы и о формировании своих затрат. Все что доступно пользователю на шаге $k+1$ – это история затрат на разных путях, соответствующих его корреспонденции, на всех предыдущих шагах $\left\{ l_p^r = \left\{ G_p\left(x^r\right) \right\}_{p \in P_w} \right\}_{r=1}^{k}$. Для простоты рассуждений мы не зашумляем эту информацию, считая что доступны точные значения имевших место затрат. Все последующие рассуждения (в виду общности выбранной в разделе 2) можно обобщить и на случай зашумленных данных (детали мы вынуждены здесь опустить).

Допуская, что $0 \leq \left\{ l_p^k \right\} \leq \breve{M}$ могут выбираться враждебно, пользователь стремиться действовать оптимальным образом, то есть так, как предписывает стратегия из примера 3 (с $i = p$, $n = |P_w|$). Заметим, что при некоторых дополнительных оговорках (см. п. 5.4 [1]) случайный выбор пути (согласно примеру 3) может быть осуществлен за время $\mathrm{O}(|E|)$, что не зависит от $n$, которое может быть намного больше (например, для манхетенской сети, см. п. 5.4 [1]).

Представим себе, что остальные пользователи ведут себя аналогичным образом, но независимо (в вероятностном плане) друг от друга. Тогда в пределе $\bar{N} \to \infty$ такая стохастическая марковская динамика в дискретном времени вырождается в детерминированную динамику в дискретном времени [16], описываемую итерационным процессом из примера 2 с $j = w$, $z^j = \left\{ x_p \right\}_{p \in P_w}$, $n_j = |P_w|$, $M = \breve{M}\sqrt{m}$,

$$\alpha_j = M^{-1}\sqrt{2\ln n_j / N}, \ N\mathrm{Regret}_N \leq \max_{j=1,\ldots,m} \alpha_j^{-1} \left( R^2 + \frac{1}{2} M^2 N \max_{j=1,\ldots,m} \alpha_j^2 \right),$$

для задачи не онлайн оптимизации:

$$\Psi(x) \to \min_{x \in X}, \tag{4}$$

$$f_k(x) \equiv \Psi(x), \ \bar{x}^N = \frac{1}{N} \sum_{k=1}^{N} x^k, \ \Psi_* = \Psi(x_*),$$

$$\Psi\left(\bar{x}^N\right) - \Psi_* \leq \mathrm{Regret}_N \leq \frac{M}{\sqrt{N}} \frac{\max_{j=1,\ldots,m} \left\{ \ln n_j \right\}}{\sqrt{2 \min_{j=1,\ldots,m} \left\{ \ln n_j \right\}}} \left( \sum_{j=1}^{m} d_j^2 + 1 \right).$$

Решение задачи (4) иногда называют равновесием Нэша(–Вардропа) в описанной популяционной игре [4, 17], соответствующей модели Бэкмана равновесного распределения потоков по путям [2]. Для простоты формулировок, будем далее считать, что решение единственно.

Введем теперь схожий процесс (совпадающий с описанным ранее в пределе $\bar{N} \to \infty$): дискретный аналог имитационной логит динамики с произвольным параметром $\alpha > 0$, популярной в эволюционной теории игр [4]. Пусть отрезок времени $[0,T]$ разбит на $T\tilde{N} \gg 1$ одинаковых отрезков, соответствующих шагам. На каждом шаге $k = 1,\ldots,T\tilde{N}$ каждый пользователь корреспонденции $j = w \in OD$ независимо от всех остальных пользователей с вероятностью $\tilde{N}^{-1}$ принимает решение выбрать потенциально новую стратегию (маршрут следования) согласно распределению вероятностей $i \in P_j$ (в действительности, тут требуются некоторые оговорки на случай когда, $x_i^k = 0$, мы опускаем здесь эти детали, за подробностями отсылаем к монографии [4])

$$p_i^{k+1} = \frac{x_i^k \exp\left(-\alpha G_i\left(x^k\right)\right)}{\sum_{l \in P_j} x_l^k \exp\left(-\alpha G_l\left(x^k\right)\right)},$$

а с вероятностью $1 - \tilde{N}^{-1}$ – использовать стратегию предыдущего шага. Аналогично действуют пользователи, принадлежащие другим корреспонденциям $j = 1,...,m$. Тогда в пределе $\tilde{N} \to \infty$ эта динамика превратится на отрезке $[0,T]$ в имитационную логит динамику в непрерывном времени [4, 16], в которой с каждым пользователем связан свой (независимый) Пуассоновский процесс с интенсивностью 1. В моменты скачков процесса пользователь принимает решение о потенциальной смене маршрута следования согласно распределению вероятностей $i \in P_j$, $j = 1,...,m$

$$p_i(t) = \frac{x_i(t) \exp\left(-\alpha G_i\left(x(t)\right)\right)}{\sum_{l \in P_j} x_l(t) \exp\left(-\alpha G_l\left(x(t)\right)\right)}.$$

При $T \to \infty$ описанный эргодический марковский процесс выходит на стационарную вероятностную меру [4]

$$\sim \exp\left(-\bar{N} \cdot \left(\Psi(x) + o(1)\right)\right),$$

которая при $\bar{N} \to \infty$ экспоненциально концентрируется в окрестности решения задачи (4).

Если описанные предельные переходы выполнить в обратном порядке: сначала $\bar{N} \to \infty$, потом $T \to \infty$, то марковский процесс, отвечающий имитационной логит динамике выродится в СОДУ $i \in P_j$, $j = 1,...,m$

$$\frac{dx_i(t)}{dt} = d_j \frac{x_i(t) \exp\left(-\alpha G_i\left(x(t)\right)\right)}{\sum_{l \in P_j} x_l(t) \exp\left(-\alpha G_l\left(x(t)\right)\right)} - x_i. \qquad (5)$$

Эта динамика (на внутренности инвариантного относительной этой динамики множества $X$) имеет глобальным аттрактором – неподвижную точку, определяемую как решение задачи (4). Более того, СОДУ (5) имеет функцию Ляпунова $\Psi(x)$ [4] (это общий факт: функционал Санова – является функционалом Больцмана [18]), причём [5]

$$\Psi\left(\frac{1}{T}\int_0^T x(t)\,dt\right) - \Psi_* \le \frac{1}{\alpha T}\sum_{j=1}^m d_j^2 \ln n_j.$$

Заметим также, что СОДУ (5) можно понимать, как непрерывный аналог (см., например, [19]) примера 2.

## **Литература**